\documentclass{article}
\usepackage{amssymb,latexsym}
\usepackage{amsthm}
\usepackage{psfrag}
\usepackage[all]{xy}
\newcommand{\Ric}{\ensuremath{\mathrm{Ric}}}

\newtheorem{thm}{Theorem}[section]
\newtheorem{lemma}[thm]{Lemma}
\newtheorem{dfn}[thm]{Definition}
\newtheorem{ex}[thm]{Example}
\newtheorem{cor}[thm]{Corollary}

\begin{document}

\title{Noncompact Manifolds with Nonnegative Ricci Curvature}
\author{William C. Wylie}
\date{}
\maketitle
\begin{abstract}
Let $(M,d)$ be a metric space.  For $0<r<R$, let $G(p,r,R)$ be the group obtained by considering all loops based at a point $p \in M$ whose image is contained in the closed ball of radius $r$ and identifying two loops if there is a homotopy betweeen them that is contained in the open ball of radius $R$.  In this paper we study the asymptotic behavior of the $G(p,r,R)$  groups of complete open manifolds of nonnegative Ricci curvature.  We also find relationships between the $G(p,r,R)$ groups and tangent cones at infinity of a metric space and show that any tangent cone at infinity of a complete open manifold of nonnegative Ricci curvature and small linear diameter growth is its own universal cover.  \end{abstract}

\section{Introduction}
One of the most fundamental areas of Riemannian geometry is the study of the relationship between curvature and topological structure. 

The case of open complete manifolds with nonnegative sectional curvature is well understood.  In fact, Cheeger-Gromoll's Soul Theorem states that any complete open manifold with nonnegative sectional curvature is diffeomorphic to a normal bundle over a compact, totally geodesic submanifold $S$ called the soul \cite{ChGr}.  In particular, every open manifold that supports a complete metric of nonnegative sectional curvature has finite topology.  Abresch-Gromoll \cite{AbGr}, Shen \cite{Sh}, and Sormani \cite{Sor2} have shown that certain exceptional subclasses of manifolds with nonnegative Ricci curvature also have finite topological type.  However, there are many examples of complete Riemannian manifolds with nonnegative Ricci curvature for which the Soul Theorem does not hold.  In fact, many of these examples do not even have finite topological type. See for example \cite{GrMe}, \cite{Men3},  \cite{Men2}, \cite{ShYa}, and \cite{Wei}.  

On the other hand, it is still unknown whether  every open manifold that supports a complete metric of nonnegative Ricci curvature must have a finitely generated fundamental group.  This was first conjectured to be true in \cite{Mil} by Milnor and a lot of interesting work has been done on the problem by Anderson \cite{And}, Li \cite{Li}, Sormani \cite{Sor1}, and  Wilking \cite{Wilk} among others.   In this paper we extend some results of Sormani \cite{Sor1} and Xu, Wang, and Yang \cite{YXW}  to obtain an understanding of how the fundamental group sits geometrically in the manifold in a sense that we will explain now. 

\begin{dfn}\emph{Given a point $p \in M$ and $0<r<R$ define $G(p,r,R)$  as the group obtained by taking all the loops based at $p$ contained in the closed ball of radius $r$ and identifying two loops if there is a homotopy between them that is contained in the open ball of radius $R$. We will refer to these groups as the geometric semi-local fundamental groups of $M$ at $p$.} \end{dfn}

The definition of the geometric semi-local fundamental groups is motivated by the concept of relative $\delta$-covers defined by Sormani and Wei. (\cite{SorWei2}, Defn 2.6).   It may seem strange to the reader to consider inner balls which are closed and outer ones which are open.  The definition is given as above because we can find a nice characterization of $G(p,r,R)$ as a subgroup of a group of deck transformations (See Corollary   \ref{Good} and Lemma \ref{Halfway}).  Note that the there is an obvious  map from $G(p,r,R)$ to $\pi_1(M,p)$ induced by the inclusion of $\overline{B_p(r)}$ into $M$.  In this paper $\pi_1(M) \cong G(p,r,R)$ will mean that this induced map is a group isomorphism.

Note that the geometric semi-local fundamental groups depend heavily on the metric structure of $M$ and not just the topology.  Even a simply connected manifold may have very complicated geometric semi-local fundamental  groups.   However, for  any complete manifold of nonnegative sectional curvature, there is some large number $r_0$ (that depends on the manifold and point $p$), so that all the groups $G(p, r, R)$ are isomorphic to $\pi_1(M)$ for all $r>r_0$ (see Example \ref{examp3}) .  In this paper  we show that there is a large class of manifolds with nonnegative Ricci curvature for which we have similar asymptotic control on the geometric semi-local fundamental groups.  Specifically,

\begin{thm} \label{Thm1}
There exists a constant $S_n = \displaystyle \frac{1}{4} \frac{1}{n-1} \Big( \frac{n-2}{n-1} \Big)^{n-1}$ and a function $g_n: (0, 2S_n] \longrightarrow [0,\infty]$ with $\displaystyle \lim_{s \rightarrow 0} g_n(s) = 1$ such that if $M^n$ is a complete open Riemannian manifold with nonnegative Ricci curvature and $p \in M$ with
\[ \displaystyle \limsup_{r \rightarrow \infty} \frac{D_p(r)}{r} < s \]
for some $s \in [0,2S_n]$ then there is some $R$ large enough  (depending on $M$) so that 
\[\pi_1(M) \cong G(p,r,g_n(s)r) \qquad \forall r > R.\]
\end{thm}

In Theorem \ref{Thm1} $D_p(r)$ is the ray density function, see Definition \ref{RayDensity}.  Roughly, Theorem 1 states that the geometric semi-local fundamental groups can be controlled if, near infinity, every point is sufficiently close to a ray.  A rougher idea of a space having many rays near infinity is the concept of a space being asymptotically polar (See Definition \ref{AssPole}).  In this case we have a similar theorem. 
\begin{thm} \label{Thm2}
Let M be a  complete open Riemannian manifold of nonnegative Ricci curvature which is asymptotically polar.  Then, for every $\varepsilon > 0$ and every point $p \in M$, there is  R large enough (possibly depending on p, M and $\varepsilon$) so that
\[ \pi_1(M) \cong G(p, r, (1+\varepsilon)r) \qquad \forall r>R.\]
\end{thm}

There are many examples that satisfy the hypotheses of Theorems \ref{Thm1} and \ref{Thm2}.  Specifically, Cheeger and Colding \cite{ChCo} have shown that manifolds with nonnegative Ricci curvature and Euclidean volume growth are asymptotically polar; and Sormani \cite{Sor3} has shown that manifolds with nonnegative Ricci curvature and linear volume growth satisfy the hypotheses of either theorem.  In \cite{Men1} Menguy gives an example of a manifold with nonnegative Ricci curvature which is not asymptotically polar.  However, the conclusion of Theorem \ref{Thm2} holds for this manifold.  The author is unaware of any examples of manifolds of nonnegative Ricci curvature which do not satisfy the conclusion of Theorem \ref{Thm2}.

The conclusions to Theorems \ref{Thm1} and \ref{Thm2} specifically imply that $\pi_1(M)$ is finitely generated.  Under the hypothesis of Theorem \ref{Thm1} Xu, Yang, and Wang (Theorem 3.1, \cite{YXW}) have shown that $\pi_1(M)$ is finitely generated.  As in this paper, their work is based on original work of   Sormani who proved that $\pi_1(M)$ is finitely generated when $M$ satisfies the hypothesis for a smaller constant $S_n$ (Theorem 1, \cite{Sor1}).  Similarly, Theorem \ref{Thm2} is a strengthening of Sormani's  Pole Group Theorem (Theorem 11, \cite{Sor1})  in which she proves that every asympototically polar, complete, open manifold with nonnegative Ricci curvature has finitely generated fundamental group.  

Theorems \ref{Thm1} and \ref{Thm2} are stronger than the results above because they not only control the generators of $\pi_1(M)$ but also the relations.  To see what new information these stronger conclusions give us recall the following definition. 
\begin{dfn} (\cite{Sp} p 62, 82) \emph{A path connected covering space $\tilde{Y}$ of another path connected topological space $Y$ is  the universal cover of $Y$  if $\tilde{Y}$  covers every other path connected cover of $Y$ and the covering projections form a commutative diagram.}  \end{dfn}

It is a well known fact that if $Y$ is  semi-locally simply connected then the universal cover exists and is the unique simply connected cover of $Y$ (\cite{Sp}, p 87, cor 4). However, if $Y$ is not semi-locally simply connected, the universal cover may or may not exist and will not be simply connected(\cite{SorWei1}, Section 2).   

As mentioned above, groups similar to the geometric semi-local fundamental groups  were used by Sormani and Wei in \cite{SorWei2} to study the topology of limit spaces of manifolds with Ricci curvature bounded below. Theorem \ref{Thm3} is a direct application of the work and can be thought of as a kind of partial converse to Theorem \ref{Thm2}.

\begin{thm} \label{Thm3}
Let $X$ be a complete length space and let $x \in X$.  If there exists positive numbers $k>1$ and $R$ so that 
\[ \pi_1(X)  \cong G(p,r,kr) \qquad \forall r>R\]
then any tangent cone at infinity of $X$ (see Definition \ref{TConeInfin}) is its own universal cover.\end{thm}

Theorem \ref{Thm3} is not true if we only assume $X$ has finitely generated fundamental group. Thus, the extra control gained by Theorem \ref{Thm1} does yield additional information.  We say that a manifold that satisfies the hypothesis of Theorem \ref{Thm1} for some $s \in (0,S_n]$ has small linear diameter growth with respect to ray density.   Directly applying Theorems \ref{Thm1} and \ref{Thm3} we obtain, 
\begin{cor}\label{Cor1}
If $M^n$ is an open manifold with nonnegative Ricci curvature and small linear diameter growth with respect to ray density then any tangent cone at infinity of $M$ is its own universal cover.
\end{cor} 

As mentioned above, in the case where the tangent cone at infinity is semi-locally simply connected the conclusion to Theorem \ref{Thm3} and Corollary \ref{Cor1} is equivalent to saying the tangent cone at infinity is simply connected.   However, it is unknown whether the tangent cones at infinity of a manifold with nonnegative Ricci curvature must be semi-locally simply connected.  In \cite{SorWei2} Sormani and Wei  do show that the Gromov Hausdorff limit of a sequence of spaces with a universal lower bound on Ricci curvature has a universal cover, although we do not require this result to prove Corollary \ref{Cor1}.

The paper is organized as follows.  In the next section we introduce the nullhomotopy radius function $\rho$ which measures how different the geometric semi-local fundamental groups are from $\pi_1(M)$.  We also give examples and discuss some basic properties of the geometric semi-local fundamental groups and the nullhomotopy radius.    Just as in \cite{Sor1} and \cite{YXW} the proofs of Theorems \ref{Thm1} and \ref{Thm2} are the result of two lemmas,  the Halfway Lemma and the Uniform Cut Lemma.  In Sections 3 and 4 we give proofs of versions of the two lemmas for the geometric semi-local fundamental groups.  The main difficulty here is that we work on the universal cover of an open metric ball which may not be a complete metric space.   Thus, these sections consist of a series of technical lemmas to work around this difficulty.  In Section 5 we give the proof of Theorem \ref{Thm1}.  In the final section we discuss Theorems \ref{Thm2} and \ref{Thm3}.

\textbf{Acknowledgements.}    I would like to thank Professor Guofang Wei for all of her advice, encouragement, and support and Professor Christina Sormani for many helpful and insightful suggestions in revising this draft.
 
\section{Nullhomotopy radius}

In this section we introduce the nullhomotopy radius and discuss its relationship with the geometric semi-local fundamental  groups.  We also give some basic examples of how the geometric semi-local fundamental groups and the nullhomotopy radius  behave.  The definitions can be applied to any length space $X$.   We also show that the nullhomotopy radius is finite  for complete manifolds with nonnegative Ricci curvature (Corollary \ref{MilGro}) which will be a very important fact in Sections 5 and 6.  To end the section we prove Lemma \ref{Helper} which we will apply in the proofs of Theorems \ref{Thm1} and \ref{Thm2}.

To motivate the definition of nullhomotopy radius let us consider two loops $\gamma_1$ and $\gamma_2$ based at a point $p \in X$ with the image of both loops contained in $\overline{B_p(r)}$.  We would like to know whether the two loops are homotopic in $X$.  This is equivalent to asking whether the loop $\omega =\gamma_1 * \gamma_2^{-1}$ is nullhomotopic.  Given $R>r$, it may be that $\omega$ is nullhomotopic but there is no nullhomotopy  of $\omega$ contained in $B_p(R)$ even for an $R$ much larger than $r$. We would like to measure how much bigger than $r$ we need to make $R$ in order to check that $\gamma_1$ and $\gamma_2$ are homotopic.

\begin{dfn}\label{def1}
\emph{Let $\Omega_{p,r}$ be the set of all $R \in \mathbb{R}$ so that all loops in $\overline{B_p(r)}$ that are nullhomotopic in $X$ are also nullhomotopic in $B_p(R)$.} \end{dfn}

If $R \in \Omega_{p,r}$ then to check whether a loop in $\overline{B_p(r)}$ is nullhomotopic we only need to check for homotopies that are contained in $B_p(R)$.  The nullhomotopy radius will  be the smallest such $R$. 

\begin{dfn} \emph{The nullhomotopy radius function at $p \in X$ is the function $\rho_p: \mathbb{R}^+ \rightarrow \mathbb{R}^+ \cup \{\infty\}$ defined as}
\[\rho_p(r) =  \left\{
\begin{array}{ll}
 \inf \{R \in \Omega_{p,r}\} & \mbox{if $\Omega_{p,r} \neq \emptyset$} \\
\infty & \mbox{if $\Omega_{p,r} = \emptyset$}\\
\end{array}
\right. .\]
\emph{We say that $\rho_p(r)$ is the nullhomotopy radius at $p$ with respect to $r$.}
\end{dfn}
An equivalent definition of $\Omega_{p,r}$ is as the set of $R$ so that the natural map from $G(p,r,R)$ to $Im \{ \pi_1(\overline{B_p(r)}) \longrightarrow \pi_1(M) \}$ is an isomorphism. Thus, $\rho_p(r)$ is the smallest number so that $G(p,r,R)$ is isomorphic to $Im \{ \pi_1(\overline{B_p(r)}) \longrightarrow \pi_1(M) \}$ for all $R > \rho_p(r).$  When it is clear which basepoint we are using  we will suppress the point and write $\rho(r)$.  Let us illustrate the behavior of these quantities with a few basic examples.

\begin{ex} \emph{Let $S^1_{\varepsilon}$ be the round circle of circumference $\varepsilon$.  Let N be a flat $S^1_{\varepsilon} \times [0, \infty)$ with the upper hemisphere of a round $S^2_{\varepsilon}$  glued to $S^1_\varepsilon \times \{0\}$ and let $p\in N$ such that the distance, $l$, from the north pole of the $S^2$ is larger than $\frac{\varepsilon}{2}$ .   $N$ is simply connected.   However, there is a geodesic loop of length $\varepsilon$ at $p$ that is not nullhomotopic inside any ball of radius less than $l$.  Namely the curve that wraps around the cylinder is only nullhomotopic via homotoping it  over the top of the attached sphere.  Thus  the nullhomotopy radius at $p$ with respect to $r$ is}
\[\rho(r) =  \left\{
\begin{array}{ll}
r& \mbox{if $r<\frac{\varepsilon}{2}$} \\
l & \mbox{if $\frac{\varepsilon}{2} \leq r < l$}\\
r & \mbox{if $r \geq l$}\\
\end{array}
\right. .\]
\end{ex}

\begin{ex}\emph{ For another simple example, consider $N = T^n \times \mathbb{R}$ where $T^n$ is the n-torus $T^n = S^1_\varepsilon \times S^1_\varepsilon \times \cdots \times S^1_\varepsilon$.  Then the fundamental group of $N$ is the free abelian group on n generators.  However, $G\left(p, \frac{\varepsilon}{2}, R\right)$ is the free group on n generators for any $R< \frac{\sqrt{n}}{2} \varepsilon$.  This is because the commutator of the minimal geodesic loops that represent the generators of $\pi_1(N)$ are contained in the ball of radius $\frac{\varepsilon}{2}$ and are nullhomotopic in $N$ but are not nullhomotopic in any ball of radius less than $ \frac{\sqrt{n}}{2} \varepsilon$.    Thus the nullhomotopy radius at $p$ with respect to $r$  at $p$ is}
\[\rho(r) =  \left\{
\begin{array}{ll}
r& \mbox{if $r<\frac{\varepsilon}{2}$} \\
 \frac{\sqrt{n}}{2} \varepsilon & \mbox{if $\frac{\varepsilon}{2} \leq r <  \frac{\sqrt{n}}{2} \varepsilon$}\\
r & \mbox{if $r \geq \frac{\sqrt{n}}{2} \varepsilon$}\\
\end{array}
\right. .\]
\end{ex}

\begin{ex}\label{examp3}\emph{More generally, let $N$ be any complete noncompact Riemannian manifold with nonnegative sectional curvature.  Then, by a theorem of Sharafutdinov \cite{Shar} (see Thm 2.3 in \cite{Yim}) , $N$ deformation retracts onto a compact soul $S$ and this deformation retraction is distance nonincreasing.  For any $p$ in $M$  take a large enough $r$ such that $B_p(r)$ contains $S$.  By following Sharafutdinov's deformation retraction, any loop in $B_p(r)$ can be homotoped into $S$ while staying inside $B_p(r)$.     Thus,  if $M$ has nonnegative sectional curvature, there always exists a large enough $R$ so that $\rho(r) = r$ for all $r>R$.  The first two examples show that this R may be very large and depend upon the point we choose.}
\end{ex}

\begin{ex} \emph{We can also construct simple spaces where the nullhomotopy radius is not well behaved.  Consider the standard, flat $xy$-plane sitting in $\mathbb{R}^3$ with standard Euclidean coordinates.  For each positive integer $n$, remove a small disc in the $xy$-plane around each point  $(n,0,0)$  and glue in its place a long capped, flat cylinder with the cap of the cylinder at the point $(n,0,10^n)$ .  Let $M$ be the resulting simply connected metric space and let $M$ be the point $(0,0,0)$.  Then for each ball of radius $n + 1/2$ around $p$ in $M$ the loop that  wraps around the glued in capped cylinder  is nullhomotopic in $M$ but not inside any metric ball centered at the origin of radius less that $n+10^n$.  Thus $\rho(n+1/2) > n + 10^n$ for all $n$ and $G(p,n+1/2, n+10^n)$ is not isomorphic to $\pi_1(M)$ for any $n$.}  \end{ex}

\begin{ex} \label{infinite}  \emph{Bowditch and Mess \cite{BoMe} and  Potyagailo \cite{Po}  have shown that there are examples of complete hyperbolic manifolds such that $\rho(r) = \infty$ for all $r \geq 1.$ This is also pointed out by Sormani and Wei, see Example 4.1 in \cite{SorWei2} .}
 \end{ex}
 
In general the nullhomotopy radius may be infinite.  However, the nullhomotopy radius will be finite in the case of nonnegative Ricci curvature. To prove this we use a very deep theorem that follows from work of Gromov and Milnor.

\begin{thm} \label{Gromov, Milnor} \textbf{(Gromov, Milnor)} \cite{Gr}, \cite{Mil} If $M$ is a complete manifold with nonnegative Ricci curvature, then every finitely generated subgroup of $\pi_1(M)$ has a nilpotent subgroup of finite index.
\end{thm}

This theorem is of particular interest to us because it implies that every finitely generated subgroup of $\pi_1(M)$ is also finitely presented.  This combined with the following general observation shows that  the nullhomotopy radius is finite for nonnegative Ricci curvature.

\begin{lemma} \label{FinPres}
Let $X$ be a complete length space.  If $Im \{ \pi_1(\overline{B_p(r)}) \longrightarrow \pi_1(X) \}$ is finitely presented then $\rho(r)$ is finite. 
\end{lemma}
\begin{proof}
Suppose that $Im \{ \pi_1(\overline{B_p(r)}) \longrightarrow \pi_1(X) \}$ is finitely presented.   Note that $G(p,r,2r)$ is finitely generated.   Let $\{\gamma_1, \gamma_2, \cdots , \gamma_k\}$ be a finite collection of loops in $\overline{B_p(r)}$ such that $G(p,r,2r) = \langle \{ [\gamma_1], [\gamma_2], \cdots , [\gamma_k] \} \rangle $. Then $\{ [\gamma_1], [\gamma_2], \cdots , [\gamma_k] \}$ is a finite set of generators for $Im \{ \pi_1(\overline{B_p(r)}) \longrightarrow \pi_1(X) \}$. Since it is finitely presented, we can write a presentation for $Im \{ \pi_1(\overline{B_p(r)}) \longrightarrow \pi_1(X) \}$ of the form
\[ \left \langle [\gamma_1], [\gamma_2], \cdots , [\gamma_k] | R_1, R_2, \cdots R_l \right \rangle \]
Each $R_j$ can be represented by a homotopy involving the representative loops $\{ \gamma_1, \gamma_2, \cdots, \gamma_k \}$. For each $j$,  let $H_j$ be this homotopy.  Since there are only finitely many homotopies we know that there exists an $R>2r$ such that
\[ \displaystyle \bigcup^l_{j=1} Im(H_j) \subset B_p(R) \]
Let $\sigma$ be a loop contained in $\overline{B_p(r)}$ with $[\sigma] = 0$ in $\pi_1(X).$  Then, since $R> 2r$, $[\sigma] = [w]$ in $G(p,r,R)$, where $w$ is some word in the $\gamma_i$s.  $[w] = 0$ in $\pi_1(X)$ so it can be written as a product of the $R_j$s.  But  $w$ can be homotoped to the constant map via the homotopies $H_j$. Since the image of all these homotopies are contained in $B_p(R)$ we see that $[w] = 0$ in $G(p,r,R)$.  Since we have started with an arbitrary nullhomotopic loop in $\overline{B_p(r)}$, $R \in \Omega_{p,r}$.
\end{proof}
 
 \begin{cor} \label{MilGro}
If $M^n$ is a complete Riemannian manifold with nonnegative Ricci curvature, then $\rho(r) < \infty$ for all $r$.
 \end{cor}
 
There is one more basic fact about the nullhomotopy radius that we will state here.  The proofs of   Theorems \ref{Thm1} and \ref{Thm2} are applications of lemma \ref{Helper}.  
  
 \begin{lemma} \label{Helper}
Let $X$ be a complete length space.  If there exists positive numbers  $L$, $k$, and $N_0$ such that  $\rho(L) < \infty$ and $G(p,L,kr) \cong G(p,r,kr)$ for all $r>N_0$.  Then  there exists $R_0$ such that $\pi_1(M) \cong G(p,r,kr)$ for all $r>R_0$.  \end{lemma}
\begin{proof}
Take $r$ so that $r>\frac{\rho(L)}{k}$ and $r>L$.  Let $\sigma$ be a loop in $\overline{B_p(r)}$ with $[\sigma] = 0$ in $\pi_1(X)$.  Then by hypothesis, $[\sigma] = [\alpha]$ in $G(p,r,kr)$ where $\alpha$ is contained in $\overline{B_p(L)}$.  But then, since $kr > \rho(L)$, $[\alpha] = 0$ in $G(p,r,kr)$.  Thus the natural map from $G(p,r,kr)$ to $\pi_1(M)$ is one to one.  The hypothesis clearly implies that the map is onto.
\end{proof}
 
\section{The Halfway Lemma for G(p,r,R)}

In this section we establish the Halfway Lemma for $G(p, r, R)$ which is motivated by Sormani's Halfway Lemma in \cite{Sor1}.   In this section $N$ is a complete Riemannian manifold without boundary.  We do not require a curvature bound in this section. 

Fix $p \in N$.  $B_p(R)$, the open metric ball of radius $R$,  is an open subset of $N$ and is thus semi-locally simply connected. Let $\widetilde{B_p(R)}$ be the universal cover of $B_p(R)$ and fix a lift of $p$ to $\widetilde{B_p(R)}$, $\tilde{p}$ .  $G(p,r,R)$ is a subgroup of $\pi_1(B_p(R))$   thus,  we can identify $G(p,r,R)$ with a subgroup of the deck transformations of  $\widetilde{B_p(R)}$ in the standard way.    Note that $\overline{B_p(R)}$ may not be semi-locally simply connected, this is the reason for using the open outer ball in the definition of the geometric semi-local fundamental groups.  We also equip $\widetilde{B_p(R)}$ with the covering metric coming from $B_p(R)$.   $\widetilde{B_p(R)}$  is then a Riemmanian manifold without boundary which is not complete.  Thus, large closed metric balls in $\widetilde{B_p(R)}$ may not be compact and we do not know, apriori, that there are only a finite number of deck transformations that  move the basepoint a given distance.  However, we can argue that this is true for small enough distances.
  
\begin{lemma} \label{Finite}For any $0<r<R$ there exists a $\delta_0$ such that the set $\{g \in \pi_1\left(B_p\left(R\right), p\right) : d(\tilde{p}, g\tilde{p}) \leq 2r + \delta_0 \}$ is finite.
\end{lemma}
\begin{proof}
Let $\delta_0 > 0$ such that $r+\delta_0 <R$ and let $D$ be a smooth compact region in $M$ such that $B_p(r+\delta_0) \subset D \subset B_p(R)$.  Fix $\tilde{p}'$ as a lift of $p$ to $\tilde{D}$.  Let $i_*: \pi_1(D,p) \longrightarrow \pi_1\left(B_p\left(R\right), p\right)$ be the induced map coming from the inclusion.  Then $i_*$ maps the set $\{ h \in \pi_1(D, p) : d(\tilde{p}', h\tilde{p}') \leq 2r + \delta_0 \}$ onto the set $\{g \in \pi_1\left(B_p\left(R\right), p\right) : d(\tilde{p}, g\tilde{p}) \leq 2r + \delta_0 \}$.  Since $D$ is a compact length space so is $\tilde{D}$ and thus $B_{\tilde{p}'}(2r + \delta_0)$ is compact and  the set $\{ h \in \pi_1(D, p) : d(\tilde{p}', h\tilde{p}') \leq 2r + \delta_0 \}$ is finite.
\end{proof}

We now give the following characterization of $G(p,r,R)$ as the  subgroup of the group of deck transformations of $\widetilde{B_p(R)}$ generated by deck transformations that move the basepoint short distances.  

\begin{cor}  \label{Good} $G(p, r, R) = \langle \{ g \in \pi_1\left(B_p\left(R\right), p\right) : d(\tilde{p}, g\tilde{p}) \leq 2r \}  \rangle$ \end{cor}
\begin{proof}
Let $g \in G(p, r, R)$,  then there is a loop $\gamma: [0, L] \longrightarrow B_p(r)$ with $[\gamma] = g$ in $\pi_1(B_p(R), p)$.  Assume $\gamma$ is parametrized  by arclength. Fix $\delta > 0$ and 
 Let $0 =t_0 < t_1<t_2< \cdots <t_k=L$ such that $t_{i+1}-t_i < \delta$.  Let $\sigma_i$ be the minimal geodesic in $B_p(r)$ from $p$ to $\gamma(t_i)$ and let $\alpha_i$ be the loop based at $p$ which traverses $\sigma_i$ from $p$ to $\gamma(t_i)$ then proceeds along $\gamma$ to $\gamma(t_{i+1})$ then returns to $p$ via $\sigma_{i+1}$.  Then $L(\alpha_i) \leq 2r + \delta$ where $L(\alpha_i)$ denotes the length of $\alpha_i$.   Let  $h_i =[\alpha_i]$, then $d(\tilde{p}, h_i(\tilde{p})) \leq 2r+ \delta$.
But $g = h_1h_2 \cdots h_k$.  So we have proven that 
\[G(p, r, R) \subset \langle \lbrace g : d(\tilde{p}, g(\tilde{p})) \leq 2r +\delta \rbrace \rangle \qquad \forall \delta\]

But  by Lemma \ref{Finite} , the set $\lbrace g : d(\tilde{p}, g(\tilde{p})) \leq 2r+\delta_0 \rbrace$ is finite for a small $\delta_0$. So there is a possibly smaller $\delta_0$ so that $\lbrace g: d(\tilde{p}, g(\tilde{p})) \leq 2r+\delta_0 \rbrace = \lbrace g : d(\tilde{p}, g(\tilde{p})) \leq 2r \rbrace $. Therefore, $G(p, r, R) \subset \langle \lbrace g : d(\tilde{p},g(\tilde{p})) \leq 2r \rbrace \rangle$.  The other inclusion is trivial.
\end{proof}

Note that using the closed ball $\overline{B_p(r)}$ is necessary in the proof of Corollary \ref{Good}.  To get a halfway generating set  we would like to take minimal geodesics in $\widetilde{B_p(R)}$ and project them down.  However, we must prove that these minimal geodesics exist.

\begin{lemma} \label{mingeod}
If $g\in \pi_1(B_p(R))$ and $d(\tilde{p}, g\tilde{p}) < 2R$ then there is a minimal geodesic in $\widetilde{B_p(R)}$ from $\tilde{p}$ to $g\tilde{p}$. 
\end{lemma}
\begin{proof}
First observe that, although $\widetilde{B_p(R)}$ is not complete, it is clear that for all $g\in \pi_1\left( B_p\left(R\right), p \right)$ the exponential map at $g\tilde{p}$ is defined on $B_0(R) \subset T_{g\tilde{p}}M$.  Therefore a minimal geodesic from $\tilde{p}$ to $g(\tilde{p})$ exists for all $g$ with $d(\tilde{p}, g\tilde{p}) < R$.   Let  $R \leq d(\tilde{p}, g\tilde{p}) < 2R$ and let $2D = d(\tilde{p}, g\tilde{p})$. $\partial{B_{\tilde{p}}(D)}$ is compact since $D<R$. Let $q$ be the point that minimizes the function $x \longrightarrow d(x, g\tilde{p})$ for $ x\in \partial{B_{\tilde{p}}(D)}$.  Then, for any $\varepsilon > 0$, there is a curve $\sigma_{\varepsilon}$ from $\tilde{p}$ to $g\tilde{p}$ with length less than or equal to $2D + \varepsilon$ and there is $t_\varepsilon \geq D$ such that $\sigma_{\varepsilon}(t_{\varepsilon}) \in \partial B_{\tilde{p}}(D)$.  Then $d(\sigma_{\varepsilon}(t_{\varepsilon}), g\tilde{p}) \leq D+\varepsilon$.   So, by the definition of $q$,  $d(q, g\tilde{p}) \leq D$.  Since $D < R$ there is a minimal geodesic from $q$ to $g\tilde{p}$ call this minimal geodesic $\gamma_2$.  Let $\gamma_1$ be a minimal geodesic from $\tilde{p}$ to $q$.  Then the curve that transverses $\gamma_1$ and then $\gamma_2$ has length less than or equal to $2D$ and therefore is a minimal geodesic from $\tilde{p}$ to $g\tilde{p}$ \end{proof}

We are now ready to prove the Halfway Lemma.  Let us first review the definition of a set of halfway generators.

\begin{dfn} \emph{A set of generators $\{g_1, g_2, ..., g_m\}$ of a group $G$ is an ordered set of generators if each $g_i$ can not be written as a word in the previous generators and their inverses.}
\end{dfn}

\begin{dfn}  \emph{Given $g \in G(p,r,R)$ we say $\gamma$ is a minimal representative geodesic loop of $g$ if $g=[\gamma]$ and $L(\gamma) = d_{\widetilde{B_p(R)}}(\tilde{p}, g\tilde{p})$.}
\end{dfn}

We can now state the Halfway Lemma for $G(p,r,R)$.

\begin{lemma}{\textbf{(Halfway Lemma)}} \label{Halfway} Let  $N$ be a complete Riemannian manifold and $p \in N$.  Then for any $0<r<R$ there exists an ordered set of generators $\{g_1, ..., g_m\}$ of $G(p, r, R)$ with minimal representative geodesic loops, $\gamma_i$ of length $d_i$ such that \[ d\left(\gamma_i(0), \gamma_i(d_i/2)\right) = d_i/2.\]
\end{lemma}

\begin{dfn}\emph{As in \cite{Sor1} We call an ordered set of generators satisfying the conclusion to Lemma \ref{Halfway} a set of halfway generators. }\end{dfn}

\begin{proof}
From the previous lemmas we know that  $G(p, r, R) = \langle \lbrace \tilde{\gamma} : d(\tilde{p}, \tilde{\gamma}(\tilde{p}) \leq 2r \rbrace \rangle$ with the generating set on the right hand side having only finitely many elements.  Therefore, we  can take $g_1$ so that $d(g_1(\tilde{p}), \tilde{p})$ is minimal among all $g_1 \in G(p, r, R)$.  Then, if $\lbrace g_1 \rbrace$ is not a generating set  of $G(p,r,R)$, consider $G(p, r, R) \setminus \langle g_1 \rangle$.  Take $g_2$ to minimize $d(g_2(\tilde{p}), \tilde{p})$ among all $g_2 \in G(p, r, R) \setminus \langle g_1 \rangle$.  By the above, $d(\tilde{p}, g_2(\tilde{p})) \leq 2r$.   Inductively we define a generating sequence with the properties that for each i, $g_i$ minimizes $d(g_i(\tilde{p}), \tilde{p})$ among all $g_i \in G(p, r, R) \setminus \langle g_1, g_2, \cdots, g_{i-1} \rangle$ and  $d_i \leq 2r$ for all $i$ where $d_i = d(\tilde{p}, g_i(\tilde{p}))$. 

 By Lemma \ref{mingeod} there is a minimal geodesic joining $\tilde{p}$ and $g_i(\tilde{p})$ in $\widetilde{B_p(R)}$.  Let $\tilde{\gamma_i}$ be this minimal geodesic and let $\gamma_i$ be the projection of $\tilde{\gamma_i}$ down to $B_p(R)$.  By our construction of the $g_i$s they have the property that if $h\in G(p, r, R)$ and $d(h(\tilde{p}), \tilde{p}) < d(g_i(\tilde{p}), \tilde{p})$ then $h \in \langle g_1, g_2, \cdots, g_{i-1} \rangle$.    To finish the lemma we need to show that $d(\gamma_i(0), \gamma_i(d_i/2)) = d_i/2$.  To do so, suppose that there is $i$ so that $d(p, \gamma_i(d_i/2)) < d_i/2$.  Let $\sigma$ be the minimal geodesic in $M$ from $p$ to $\gamma_i(d_i/2)$.  Let $h_1$ be the element of $G(p, r, R)$ represented by transversing $\gamma_i$ from $0$ to $\frac{d_i}{2}$ then following $\sigma$ back to $p$ and let $h_2$ be the element of $G(p,r,R)$ represented by transversing first $\sigma$ then $\gamma_i$ from $\frac{d_i}{2}$ to $d_i$.   Then it is clear that $d(h_1(\tilde{p}), \tilde{p})$ and $d(h_2(\tilde{p}), \tilde{p})$ are both less than $d(g_i(\tilde{p}), \tilde{p})$  which implies that $h_1$ and $h_2$ are in $\langle g_1, g_2, \cdots, g_{i-1}\rangle$.  But  $h_1h_2 = g_i$ so that $g_i \in \langle g_1, g_2, \cdots, g_{i-1}\rangle$ which is a contradiction to our choice of $g_i$.
 \end{proof}

\section{Localized Uniform Cut Lemma}

In this section $M^n$ is a complete Riemannian manifold with nonnegative Ricci curvature and dimension at least 3.  In dimension less than three, nonnegative Ricci curvature is equivalent to nonnegative sectional curvature, therefore Theorems \ref{Thm1} and \ref{Thm2} are already known to be true is this case. We wish to obtain a Uniform Cut Lemma for $G(p,r,R)$ that is similar to the ones of Sormani \cite{Sor1} and Xu, Yang and Wang \cite{YXW}.  To do this we need to apply the excess estimate of Abresch and Gromoll \cite{AbGr} to the universal cover of an open ball, which is not a complete manifold.  However, this turns out not to be a problem if we examine the proof in \cite{AbGr}.  

Recall that, given two points $p$ and $q$  in $M$, the excess function $e_{p,q}$ on $M$ is defined as 
 \[e_{p,q}(x) = d_{M}(p, x) + d_M(q, x) - d_{M}(p,q).\]
 Let $\gamma$ be a minimal geodesic from $p$ to $q$ and let $l(x) = d(x, \gamma )$.  Then the excess estimate states that 
 \begin{thm}\textbf{[Abresch-Gromoll]} \cite{AbGr}  \label{ExcEst.}
Let $M^n$ be an open complete Riemannian manifold with  $\Ric_M \geq 0$, $n \geq 3$, let $y, p, q \in M$,  $l = l(y)$ and $r_0 = d(p,y)$, $r_1 = d(q,y)$.  Assume $l < min \{ r_0, r_1 \}$.   then
\[e_{p,q}(y) \leq 2 \left( \frac{n-1}{n-2}\right) \left( \frac{1}{2} C_3 l^n\right) ^{1/(n-1)}\]
where $C_3 = \frac{n-1}{n} \left( \frac{1}{r_0 - l(y)} + \frac{1}{r_1 - l(y)} \right)$
\end{thm}
If we examine the proof of the above theorem, one notices that the key is applying the Laplacian comparison first on the small ball $B_y(l(y))$, then on the balls $B_p(R_0)$ and $B_q(R_1)$ where $R_0 = r_0 + l$ and $R_1=r_1+l$.  The Laplacian comparison is true for any compact closed ball.  Thus,  we can localize Abresch and Gromoll's proof to get the following.
\begin{thm} \textbf{[Abresch-Gromoll]} \cite{AbGr} \label{NonCompleteExcEst.}
Let $M^n$ be a  Riemannian manifold with  $\Ric_M \geq 0$, $n \geq 3$, let $y, p, q \in M$, $l = l(y)$ and $r_0 = d(p,y)$, $r_1 = d(q,y)$.  Assume $l < min \{ r_0, r_1 \}$.   If $\overline{B_p(R_0)}$ and $\overline{B_q(R_1)}$ are compact for $R_0 = r_0 + l$ and $R_1 = r_1 + l$ then 
\[e_{p,q}(y) \leq 2 \left( \frac{n-1}{n-2}\right) \left( \frac{1}{2} C_3 l^n\right) ^{1/(n-1)}\]
where $C_3 = \frac{n-1}{n} \left( \frac{1}{r_0 - l(y)} + \frac{1}{r_1 - l(y)} \right)$
\end{thm}
We now get a noncomplete version of  an important lemma of Xu, Yang, and Wang's (\cite{YXW}, Lemma 2.2).
\begin{thm}\label{NonCompleteCutLemma}
Let $M^n$ be a Riemannian manifold with $\Ric_M \geq 0$ and $n \geq 3$.  Let $\gamma$ be a minimal geodesic parametrized by arclength and L$(\gamma) = D$.  If $x \in M$ is a point with $d\left(x, \gamma(0)\right) \geq \left(\frac{1}{2} + \varepsilon_1\right)D$, and $d\left(x, \gamma(D)\right) \geq (\frac{1}{2} + \varepsilon_2)D$ and if  the closure of the balls $B_{\gamma(0)}\left(\left(\frac{1}{2} + 2\alpha \left(\varepsilon_1, \varepsilon_2\right)\right)D\right)$ and  $B_{\gamma(D)}\left(\left(\frac{1}{2}+2\alpha \left(\varepsilon_1, \varepsilon_2\right)\right)D\right)$ are compact then $d(x, \gamma(\frac{D}{2})) > \alpha(\varepsilon_1, \varepsilon_2)D$ where 
\[\alpha(\varepsilon_1, \varepsilon_2) =  \alpha_n(\varepsilon_1, \varepsilon_2) = min \left\{ \frac{1}{4},\left( \frac{\varepsilon_1 + \varepsilon_2}{2} \right) ^{\frac{n-1}{n}} \left( \frac{1}{4} \frac{n}{n-1} \left( \frac{n-2}{n-1} \right) ^{n-1} \right)  ^{\frac{1}{n}} \right\} \]
\end{thm}
\begin{proof}
Suppose on the contrary that $d(x, \gamma(D/2)) \leq \alpha(\varepsilon_1, \varepsilon_2)D$.

Then
\[l(x) \leq d(x, \gamma(D/2)) \leq \alpha(\varepsilon_1, \varepsilon_2)D \leq min\{r_0, r_1\}\]
And also
\[r_0 \leq d(\gamma(0), \gamma(D/2)) + d(\gamma(D/2), x) \leq (\frac{1}{2} + \alpha(\varepsilon_1, \varepsilon_2))D\]
and similarly
\[r_1 \leq (\frac{1}{2} + \alpha(\varepsilon_1, \varepsilon_2))D\]
By assumption, then, the closure of the balls $B_p(r_0 + l)$ and $B_p(r_1+l)$ are compact so, by Theorem \ref{NonCompleteCutLemma},
\begin{equation} \label{eqn2}
e_{\gamma(0), \gamma(D)}(x) \leq  2 \left( \frac{n-1}{n-2}\right) \left( \frac{1}{2} C_3 l^n\right) ^{1/(n-1)}
\end{equation}
On the other hand, by hypothesis,
\begin{equation}\label{eqn1}
e_{\gamma(0), \gamma(D)}(x) \geq (\varepsilon_1 + \varepsilon_2)D
\end{equation}
Note that $r_0 - l(x) > D/4$ and $r_1 - l(x) > D/4$.
Thus \begin{equation} \label{eqn3} C_3 < \displaystyle \frac{8(n-1)}{nD} \end{equation}
combining (\ref{eqn2}), (\ref{eqn1}), and (\ref{eqn3}) along with the fact that $l \leq \alpha(\varepsilon_1, \varepsilon_2)D$ we see that 
\[ (\varepsilon_1 + \varepsilon_2)D < 2 \left( \frac{n-1}{n-2}\right) \left( \frac{1}{2} \frac{8(n-1)}{nD} \left(\alpha(\varepsilon_1, \varepsilon_2)D\right)^n\right) ^{1/(n-1)} \]
Solving for $\alpha(\varepsilon_1, \varepsilon_2)$ implies that
\[\alpha(\varepsilon_1, \varepsilon_2) > \left( \frac{\varepsilon_1 + \varepsilon_2}{2} \right) ^{\frac{n-1}{n}} \left( \frac{1}{4} \frac{n}{n-1} \left( \frac{n-2}{n-1} \right) ^{n-1} \right)  ^{\frac{1}{n}} \]
Which is a contradiction to the definition of $\alpha(\varepsilon_1, \varepsilon_2)$. 

\end{proof}

We are now ready to state and prove the Uniform Cut Lemma for $G(p,r,R)$.  We set $\alpha(\varepsilon) = \alpha(\varepsilon, \varepsilon)$.  Lemma \ref{LocalUniformCutLemma} is an improvement upon a similar localized uniform cut lemma of Sormani and Wei (\cite{SorWei2}, Lemma 3.14).

 \begin{lemma} \textbf{(Localized Uniform Cut Lemma)} \label{LocalUniformCutLemma}
Let $M^n$ be a complete manifold with nonnegative Ricci curvature and $n \geq 3$ and let $\varepsilon, D, R >0$ with $R > (\frac{1}{2} + 2\alpha(\varepsilon))D$ and let $\gamma$ be a loop based at $p$ with length $L(\gamma) = D$ such that the following two conditions hold.
\begin{enumerate}
\item $L(\sigma) \geq D$ for all $\sigma$ a loop based at p such that $[\sigma] = [\gamma]$ in $G(p,r,R)$.
\item $\gamma$ is a minimal geodesic on $[0, D/2]$ and $[D/2, D]$.
\end{enumerate}
Then for any $x\in \partial B_p(TD)$ where $T \geq (\frac{1}{2} + \varepsilon)D$ 
\[d_M(x, \gamma(D/2)) \geq (T -\frac{1}{2})D + (\alpha(\varepsilon) - \varepsilon)D \]
\end{lemma}
\begin{proof}
We only need to show the above inequality for $x \in \partial B_p((\frac{1}{2}+\varepsilon)D)$.  The general one above then follows by the triangle inequality. Suppose that $d_{M}(x, \gamma(D/2)) < \alpha(\varepsilon)D$.
Let $C:[0,1] \longrightarrow M$ be the shortest path from $\gamma(D/2)$ to $x$.  Lift $\gamma$ to $\tilde{\gamma}$ in $\widetilde{B_p(R)}$. Note that since $L(C) < \alpha(\varepsilon)$ the image of $C$ is contained in $B_p((\frac{1}{2}+2\alpha(\varepsilon))D)$ so we can lift $C$ to $\tilde{C}$, a curve from $\tilde{\gamma}(D/2)$ to a point $\tilde{x}$.  Then the point $\tilde{x}$ and the minimal geodesic $\tilde{\gamma}$ will contradict Theorem \ref{NonCompleteCutLemma}.
\end{proof}

\section{Ray density and nullhomotopy radius}

In this section we present the proof of Theorem \ref{Thm1}.  Let us first review the definition of the ray density function.
\begin{dfn} \label{RayDensity} \cite{Ch} \emph{The ray density function at p is}
\[D_p(r) = \displaystyle \sup_{x \in \partial B_p(r)} \inf_{rays \gamma, \gamma(0)=p} d(x,\gamma(r)).\]
\end{dfn}

Note that the ray density is always less than or equal to the extrinsic diameter of  $\partial B_p(r)$ in $M$.  In particular, $D_p(r) \leq 2r$.

Consider a loop $\gamma$ that satisfies the hypotheses of Lemma \ref{LocalUniformCutLemma}.  Let $\sigma$ be a ray at p.  Then, by the triangle inequality and Lemma \ref{LocalUniformCutLemma}, we see that $d_M(\gamma(D/2), \sigma(D/2)) > (\alpha(\varepsilon) - \varepsilon)D$.  Let $f(\varepsilon) = \alpha(\varepsilon)-\varepsilon$.  By elementary calculus we can see that the maximum of $f$ occurs at $\varepsilon_0=(\frac{n-1}{n})^n(\frac{1}{4}\frac{n}{n-1}(\frac{n-2}{n-1})^{n-1})$, that $f(0)=0$ , that $f$ is an increasing function from $0$ to $\varepsilon_0$, and that the maximum value of $f$ is $f(\varepsilon_0)=S_n= \frac{1}{4}\frac{1}{n-1}(\frac{n-2}{n-1})^{n-1}$.

\begin{proof}[Proof of Theorem \ref{Thm1}]
To prove the theorem we need to show that the hypotheses of Lemma \ref{Helper} are satisfied.  Since Lemma \ref{MilGro} tells us that $\rho$ is  finite, we only need to show that there exists a function $g$ and an $L$ so that $G(p,L,g(s)r) = G(p,r,g(s)r)$ for all $r>L$.   To do this, set $L$ equal to the number such that $D_p(r) < rs$ for all $r>L$.  Let $\varepsilon' < \varepsilon_0$ such that $2 f(\varepsilon') = s$ and let $g(s) = 1+4(\alpha(\varepsilon'))$.  Fix $r>L$ and let $R > g(s)r$.  Then, by the Halfway Lemma, the group $G(p,r,R)$ has a set of halfway generators $\{g_1, g_2, \cdots , g_k \}$ with representative halfway loops $\gamma_i$ with $L(\gamma_i) = d_i$ and $d_i \leq 2r$.  But then since $R > g(s)r> (\frac{1}{2} + 2\alpha(\varepsilon))d_i$  we can apply the Localized Uniform Cut Lemma as in the discussion above and see that 
\[ d_M(\gamma_i(d_i/2), \sigma(d_i/2)) > (f(\varepsilon'))d_i. \]
where $\sigma$ is any ray based at $p$.  This implies that 
\[ \displaystyle \frac{D(d_i/2)}{d_i/2} \geq 2 f(\varepsilon')=s. \]
Then, by our hypothesis,  $d_i/2 \leq L$ for all i. Since $\{g_1, g_2, \cdots , g_k \}$  generate $G(p,r,g(s)r)$ this implies that $G(p,L,g(s)r) = G(p,r,g(s)r)$.
\end{proof}

\section{Tangent cones at infinity and nullhomotopy radius}
In this section we present the proofs of Theorems \ref{Thm2} and \ref{Thm3}.  We begin with the definitions of a tangent cone at infinity and an asymptotically polar manifold. 

\begin{dfn} \label{TConeInfin} \cite{GrB} \emph{Let $X$ be a complete length space.  A pointed length space, $(Y, d_{Y}, y_0)$ is called a tangent cone at infinity of $M$ if there exists $p \in M$ and a sequence of positive real numbers $\{r_i\}$ diverging to infinity such that the sequence of pointed metric spaces $(X, \frac{d_X}{r_i}, p)$ converges in the pointed Gromov-Hausdorff sense to $(Y, d_{Y}, y_0)$.} 
\end{dfn}

Note that tangent cones at infinity may not be unique and may not be metric cones.  However, by Gromov's Compactness Theorem \cite{GrB}, they do always exist for complete manifolds of nonnegative Ricci curvature.  

\begin{dfn} \cite{ChCo} \emph{A length space Y, has a pole at a point $y \in Y$ if for all $x \neq y$ there is a curve $\gamma: [0,\infty) \rightarrow Y$ such that $d_Y(\gamma(s), \gamma(t)) = |s-t|$  $\forall s,t \in [0, \infty)$ and $\gamma(0) = y$, $\gamma(t_0) = x$ for some $t_0$. }
\end{dfn}

\begin{dfn} \label{AssPole} \cite{ChCo} \emph{A complete manifold of nonnegative Ricci curvature is asymptotically polar if all of its tangent cones at infinity $(Y, d_{Y}, y_0)$ have a pole at its basepoint $y_0$.}
\end{dfn}

\begin{proof}[Proof of Theorem \ref{Thm2}]
The argument is by contradiction, that is, suppose there is a $p \in M$ and $\varepsilon_0 >0$ such that $\pi_1(M) \neq G(p, r_i, (1+\varepsilon_0)r_i)$ for some sequence $r_i \rightarrow \infty$. Then, since $\rho(L)$ is finite for all $L$,  Lemma 2.11 implies that there exists some possibly different sequence of  $\{r_i\}$ diverging to infinity such that for every $L$ there are infinitely many $r_i$ with $G(p,L,(1+\varepsilon_0)r_i) \neq G(p,r_i,(1+\varepsilon_0)r_i)$.   Let $\Gamma_i$ be the set of halfway generators of $G(p,r_i, (1+\varepsilon_0)r_i)$ and $\Gamma = \bigcup_{i=1}^{\infty} \Gamma_i$.  If there existed $L_0$ such that $L(\gamma) \leq L_0$ $\forall \gamma \in \Gamma$ then, by Lemma \ref{Halfway}, $G(p, L_0/2, (1+\varepsilon_0)r_i) = G(p, r_i , (1+\varepsilon_0)r_i)$ for all $i$, which would be a contradiction.   Therefore, there is no such $L_0$ and we can choose $\gamma_i \in \Gamma_i$ with  length  $d_i$ such that $d_i \rightarrow \infty$ as $i \rightarrow \infty$.  Also note that  $d_i \leq 2r_i$, so  \[ (1/2 + \varepsilon_0/2) d_i \leq (1+\varepsilon_0)r_i. \]   Set $\alpha(\varepsilon) = \alpha(\varepsilon, \varepsilon)$ where $\alpha$ is as in Theorem \ref{NonCompleteCutLemma}.  Take $\varepsilon$ such that $2\alpha(\varepsilon) < \varepsilon_0/2$.  Then,  by Lemma \ref{LocalUniformCutLemma}, for all $x \in \partial B_p(Td_k)$ with $T = (\frac{1}{2}+\varepsilon)$, 
\begin{equation} \label{LUCL2}
d_M(x, \gamma_i(d_i/2)) \geq \left(T -\frac{1}{2}\right)\left(d_i/2\right) + \left(\alpha(\varepsilon) - \varepsilon\right)\left(d_i/2\right) \geq \frac{\alpha(\varepsilon)}{2}d_i \end{equation}
Let $M_i = (M^n, d_M/d_i, p)$. Then this sequence has a subsequence which converges to a tangent cone at infinity $Y = (Y, d_Y, y_0)$.  We will show that $Y$ is not polar and thus contradict the assumption that $M$ is asymptotically polar.  Let $d_{GH}\left(\overline{B_p(1)} \subset M_i, \overline{B_y(1)} \subset Y\right) = 2\varepsilon_i$. Then $\varepsilon_i \longrightarrow 0$ as $i$ goes to infinity and, by the definition of Gromov-Hausdorff convergence, there exists maps 
\[F_i: \overline{B_p(1)} \subset M_i \longrightarrow \overline{B_y(1)} \subset Y\]
such that
\begin{equation} \label{AlmostContinuous}
\left|d_{M_i}(x,x') - d_{Y}(F_i(x), F_i(x'))\right| < \varepsilon_i
\end{equation}
and such that for all $y \in \overline{B_{y_0}(1)}$ there is $x_y \in \overline{B_p(1)} \subset M_i$ so that 
\begin{equation} \label{AlmostOnto}
 d_Y(F_i(x_y), y) < \varepsilon_i.
 \end{equation}
Thus, $F_i(\gamma_i(d_i/2)) \in Ann_{y_0} (1/2 - \varepsilon_i, 1/2 + \varepsilon_i)$. So the sequence $\{F_i(\gamma_i(d_i/2)) \}$ has a convergent subsequence converging to some point $y' \in \partial B_{y_0}(1/2)\subset Y$.
 
 We will show that $y'$ is the point which does not have a ray going through it.  In fact we will show that for all $y \in \partial B_{y_0}(1)$, $d_Y(y',y) > 1/2 + \frac{\alpha(\varepsilon)}{2}$.  
  
 To prove this let $y \in \partial B_{y_0}(1)$, then, by (\ref{AlmostOnto}), for any $i$ there is $x_{y,i} \in \overline{B_p(1)}$ such that $d(F_i(x_{y,i}), y) < \varepsilon_i$.  Then, by (\ref{AlmostContinuous}), \[d_{M_i}(x_{y,i}, p) \geq 1-\varepsilon_i\]
 and so
 \[d_M(x_{y,i}, p) \geq (1-\varepsilon_i)d_i\]
 Let $T_i = 1 - \varepsilon_i$.  Then, since $\varepsilon_i \longrightarrow 0$, for  large enough $i$ we know that, 
 \[T_i \geq 1/2 + \varepsilon.\]  
 By (\ref{LUCL2}) 
 \begin{eqnarray*}
 d_M(x_{y,i}, \gamma_i(d_i/2)) &\geq& (T_i-1/2)d_i +\frac{\alpha(\varepsilon)}{2}d_i \\&\geq& 1/2d_i - \varepsilon_i d_i + \frac{\alpha(\varepsilon)}{2}d_i.
 \end{eqnarray*}
 But then
 \[ d_{M_{i}}(x_{y,i}, \gamma_i(d_i/2)) \geq 1/2 - \varepsilon_i + \frac{\alpha(\varepsilon)}{2}. \]
and (\ref{AlmostContinuous}) implies
 \[ d_Y(F_i(x_{y,i}), F_i(\gamma_i(d_i/2))) \geq 1/2 - 2\varepsilon_i + \frac{\alpha(\varepsilon)}{2}.\]
 Taking the limit of both sides of the equation as $i \longrightarrow \infty$ we get that
 \[d_Y(y,y') \geq 1/2 + \frac{\alpha(\varepsilon)}{2}.\]
  
 \end{proof}

We now give the background to Theorem \ref{Thm3}.  We first review the definition of $\delta$ loops and relative $\delta$ covers introduced by Sormani and Wei in \cite{SorWei1} and \cite{SorWei2}.

\begin{dfn}
\emph{A loop $\gamma$ is called a $\delta$ loop if it is of the form $\alpha * \beta * \alpha^{-1}$ where $\beta$ is a closed path lying in a ball of radius $\delta$ and $\alpha$ is a path from $p$ to $\beta(0)$.}
\end{dfn}

\begin{dfn}
\emph{$G(p, r, R, \delta)$ is the group of equivalence classes of loops in $\overline{B_p(r)}$ where $\gamma_1$ is equivalent to $\gamma_2$ if the loop $\gamma_2^{-1} * \gamma_1$ is homotopic  in $B_p(R)$ to a product of $\delta$ loops.}
\end{dfn}
Note that here we have changed the definitions in \cite{SorWei2} slightly by changing the outer ball to an open.   This is to match with the definition of geometric semi-local fundamental groups.  However, this change will  not change the proof of the following important lemma.  Theorem \ref{Thm3} is a direct application of this lemma.

\begin{lemma} \label{WeiSormani} \textbf{(Sormani-Wei)} \cite{SorWei2}
Suppose $(M_i, p_i)$ converges to $(Y,y)$ in the pointed Gromov-Hausdorff topology.  Then for any $r<R, \delta_1<\delta_2$, there exist sequences $r_i \longrightarrow r$ and $R_i \longrightarrow R$ such that $B(p_i, r_i)$ and $B(p_i, R_i)$ converge to $B_p(r)$ and $B_p(R)$ with respect to the intrinsic metrics, and a number $N$ sufficently large depending upon $r,R,\delta_2$, and $\delta_1$ such that $\forall i > N$ there is a surjective map $\Phi_i : G(p_i,r_i,R_i, \delta_1) \longrightarrow G(y,r,R,\delta_2)$.
\end{lemma}

Theorem \ref{Thm3} follows from applying Lemma \ref{WeiSormani} to a Gromov-Hausdorff convergent sequence of rescalings of a complete length space.

\begin{lemma}\label{beef}
Let $(X, d_X)$ be a complete length space and let $x \in X$. Let $\{ \alpha_i \}$ be a sequence of positive real numbers diverging to infinity such that the sequence of pointed metric spaces  $X_i = (X, d_M/\alpha_i, p)$ converges to some $Y=(Y, d_Y, y_0)$ in the pointed Gromov Hausdorff topology.  Let $G^i(p,r,R, \delta)$ denote the relative $\delta$ group corresponding to $X_i$. If there exists positive numbers $k$ and $L$ so that 
\[ \pi_1(X) \cong G(x,r,kr) \qquad \forall r>L\]
then, for any $r', \delta > 0$, there exists $N$ sufficently large depending on $r'$ and $\delta$ such that 
\[G^{i}(x, r', (k+1)r', \delta) = 0 \qquad \forall i \geq N\]
\end{lemma}
\begin{proof}
 Let $r>L$ and $\gamma$ be a loop contained in $\overline{B_x(r)}$ then since $G(x, L, kL) \cong \pi_1(M)$ there is $\sigma$ contained in $\overline{B_x(L)}$  such that $[\sigma * \gamma ^{-1}] = 0$ in $\pi_1(X)$.,  But since $G(x, r, kr) \cong \pi_1(M)$, $[\sigma] = [\gamma]$ in $G(p,r,(k+1)r)$.  Since $\gamma$ is arbitrary and $[\sigma]=0$ in $G(x,r',(k+1)r',L)$ we have that
\[ G(x,r,(k+1)r,L) = 0 \qquad \forall r>L \] 
 In terms of $M_i$ this implies that
\[ G^{i}(x, r', (k+1)r', L/\alpha_i) = 0 \qquad \forall r' > L/\alpha_i \] 
To prove the lemma take $N$ so that $\forall i > N$ $r> r_0 /\alpha_i$ and $L/\alpha_i < \delta$.
\end{proof}

\begin{proof}[Proof of Theorem \ref{Thm3}]
Let $Y = (Y, y_0, d_Y)$ be a tangent cone at infinity of $X$.  Then combining Lemmas \ref{WeiSormani} and \ref{beef} we see that 
\begin{equation}\label{AlmostDone}
G(y_0, r, (k+1)r, \delta) = 0
\end{equation}
for all $r$ and $\delta$.  

Let $\tilde{Y}$ be any path connected cover of $Y$ and let $p: \tilde{Y} \longrightarrow Y$ be the covering map.  Then the number of sheets of the cover $\tilde{Y}$ is equal to the index of $p_*(\pi_1(\tilde{Y}))$ as a subgroup $\pi_1(Y)$.  Therefore, to show that $\tilde{Y}$ is a trivial cover, we need to show that the map $p_*$ is onto.  To do this let $\gamma$ be a closed loop based at $y_0$ in $Y$.  Then $\gamma$ is contained in $\overline{B_{y_0}(r)}$ for some $r$.  Then since $\overline{B_{y_0}(r)}$ is compact and $\tilde{Y}$ is a covering space, there is $\delta_r$ such that for all $x \in \overline{B_p(r)}$ the ball $B_x(\delta_r)$ lifts isometrically to $\tilde{Y}$.  Thus, every $\delta_r$ loop in $\overline{B_p(r)}$ is contained in the image of $p_*$.  But, by (\ref{AlmostDone}), $\gamma$ is homotopic to some product of $\delta_r$ loops and $p_*$ is onto.
\end{proof}

Mathematical Keywords: Ricci curvature, noncompact manifold, fundamental group, ray density, asymptotically polar

2000 Mathematical Subject Classification: 53C20

DEPARTMENT OF MATHEMATICS, UNIVERSITY OF CALIFORNIA, SANTA BARBARA, CALIFORNIA, 93106

\emph{email address:} wwylie@math.ucsb.edu.
\end{document}